\documentclass [a4paper] {article}

\usepackage[utf8]{inputenc}

\usepackage{amssymb,amsmath, amsfonts}

\usepackage [dvips] {graphicx,color}

\usepackage{amsthm}

\usepackage{MnSymbol}

\newcommand*{\twosquig}{%
\mathrel{\vcenter{\offinterlineskip
\hbox{$\rightsquigarrow$}\vskip-.08ex\hbox{$\rightsquigarrow$}}}}

\theoremstyle{definition}
\newtheorem{df}{Definition} 
\newtheorem*{df*}{Definition}

\theoremstyle{remark}
\newtheorem{rmq}{Remark}

\author{Stéphane \textsc{Dugowson}
\footnote{Supmeca Paris + Quartz. Email: s.dugowson@gmail.com}
}

\title {A Lax Functorial Definition of Open Dynamics}

\begin{document}

\maketitle



\paragraph{Abstract.} This paper provides a rewording in the language of lax-functors of the definition of open dynamics given in our systemic theory of interactivity previously exposed in \cite{Dugowson:20160831}, itself coming from \cite{Dugowson:20150807} and \cite{Dugowson:20150809}. 

To formulate this definition,  we need first to give the definition of a \emph{closed dynamic} on a small category $\mathbf{C}$: such a closed dynamic is a lax-functor from $\mathbf{C}$ into the large $2$-category of \emph{states sets}, \emph{transitions} and \emph{constraint order} that we  denote
\footnote{
In \cite{Dugowson:20160831}, \cite{Dugowson:20150809}, \cite{Dugowson:20150807}, \cite{Dugowson:201203} and \cite{Dugowson:201112},  the same category was denoted by $\mathbf{P}$.
}
by $\mathbf{PY}$. Using the $2$-category that we denote by $\mathbf{PY}_L$ or $\mathbf{PY}^{\underrightarrow{\scriptstyle{L}}}$, we then define in the same way $L$-multi-dynamics, where $L$ refers to any non-empty set which we call the parametric set of the considered multi-dynamic. After having defined morphisms between multi-dynamics based on different parametric sets, we define open dynamics and morphisms between them. Finally, we give the definition of a realization of an open dynamic.




\section{The bicategories $\mathbf{PY}$ and $\mathbf{PY}_L$}

\subsection{Transitions}

For any sets $U$ and $V$, we call  any map $U\rightarrow\mathcal{P}(V)$ --- or, equivalently, any binary relation $U \rightarrow V$ --- a \emph{transition from $U$ to $V$}. We often write $\varphi:U\rightsquigarrow V$ to indicate that $\varphi$ is such a transition.

Denoted $\psi\odot \varphi$, the composition of transitions $\varphi:U\rightsquigarrow V$ and $\psi:V\rightsquigarrow W$  is defined in the obvious way: for all $u\in U$ 
\[
\psi\odot \varphi (u)=\bigcup_{v\in \varphi(u)}\psi(v)\subseteq W.
\]

A transition $f:U\rightsquigarrow V$ is called
\begin{itemize}
\item \emph{deterministic}  if $card(f(u))=1$ for all $u\in U$,
\item \emph{quasi-deterministic}, if $card(f(u))\leq 1$ for all $u\in U$.
\end{itemize}

A (quasi) deterministic transition  $U\rightsquigarrow V$ is often denoted as a (partial) function $U \rightarrow V$.

\subsection{The bicategory $\mathbf{PY}$}

$\mathbf{PY}$ is defined as the $2$-category (and hence the bicategory) with objects the sets, and for each couple of sets $(U,V)$ the category $\mathbf{PY}(U,V)$ being the \emph{constraint} order on the set of  {transitions} $U\rightsquigarrow V$  defined  for all $\varphi, \psi \in \mathbf{PY}(U,V)$
by
\[\varphi\leq\psi\Leftrightarrow \varphi\supseteq \psi\]
where $\varphi\supseteq \psi$ means that for all $u\in U$, $\varphi(u)\supseteq\psi(u)$. If $\varphi\leq\psi$, we say that $\psi$ is more \emph{constraining} than $\varphi$, or that $\varphi$ is \emph{laxer} than  $\psi$. Thus, there exists a $2$-cell $\varphi \rightarrow\psi$ if and only if $\varphi$ is laxer than  $\psi$.

\subsection{The bicategory $\mathbf{PY}_L$}\label{subsection PYL}

In the same way, for any non-empty set $L$ we define the $2$-category (and hence the bicategory) denoted by $\mathbf{PY}_L$ or
$\mathbf{PY}^{\underrightarrow{\scriptstyle{L}}}$ as the one with sets as $0$-cells, and  for each couple of sets $(U,V)$ the category $\mathbf{PY}^{\underrightarrow{\scriptstyle{L}}}(U,V)$ being defined by 
\[
\mathbf{PY}^{\underrightarrow{\scriptstyle{L}}}(U,V) =
(\mathbf{PY}(U,V))^L.
\]
In other words, for a given domain $U$ and a given codomain $V$, a $1$-cells in $\mathbf{PY}^{\underrightarrow{\scriptstyle{L}}}$ is an $L$-family $\varphi=(\varphi_\lambda)$ of transitions $\varphi_\lambda:U\rightsquigarrow V$. We'll sometimes write 
$\varphi: U \twosquig_L  V $ or $ U {\stackrel{\varphi}{\twosquig}}_L V $
to say that $\varphi$ is such a family.

The composition of $1$-cells is naturally defined by $(\varphi_\lambda)\odot (\psi_\lambda)=((\varphi_\lambda\odot\psi_\lambda)_\lambda)$, and $2$-cells express the order \[(\varphi\leq \psi) \Leftrightarrow (\forall \lambda\in L, \varphi_\lambda\leq \psi_\lambda).\]

\section{Categories of closed dynamics and categories of multi-dynamics}

\subsection{The category $\mathcal{CD}_\mathbf{C}$ of closed dynamics on $\mathbf{C}$}

Let $\mathbf{C}$ be a small category viewed as a ``discrete bicategory''\footnote{For any objects $S$ and $T$, $\mathbf{C}(S,T)$ is just a set, that is a discrete category. 
}, 
and $\alpha$ be a lax-functor from $\mathbf{C}$ to $\mathbf{PY}$, that we denote\footnote{Keeping the notation we used in our previous papers about dynamics.}
\[\alpha:\mathbf{C}\rightharpoondown\mathbf{PY}.\] 
 We say that $\alpha$ is \emph{disjunctive} if for all objects $S$ and $T$ in $\mathbf{C}$,
\[
S\neq T \Rightarrow S^\alpha \cap T^\alpha=\emptyset,
\]
where we use the notation $S^\alpha$ to denote the set $\alpha(S)$.
Likewise, for $d\in\overrightarrow{\mathbf{C}}$, we put $d^\alpha=\alpha(d)$.

\begin{df} [Closed dynamics on $\mathbf{C}$] We denote by $\mathcal{CD}_\mathbf{C}$ the category of \emph{closed dynamics} on $\mathbf{C}$, that is the category
\begin{itemize}
\item whose objects are disjunctive lax-functors $\alpha:\mathbf{C}\rightharpoondown \mathbf{PY}$, 
\item whose arrows --- called \emph{dynamorphisms} --- are lax-natural transformations $\delta:\alpha\looparrowright \beta$. 
\end{itemize}
\end{df}

The small category $\mathbf{C}$ is called the \emph{engine} of the dynamics $\alpha\in Ob(\mathcal{CD}_\mathbf{C})$. The arrows $(S\stackrel{d}{\rightarrow}T)\in\overrightarrow{\mathbf{C}}$ are called \emph{durations}.


\begin{rmq}
In \cite{Dugowson:20150809}, section \textbf{§\,2.1.2}, the category $\mathcal{CD}_\mathbf{C}$ was denoted $\mathbf{DySC_{(C)}}$. Indeed, by definition of lax-functors and lax-natural transformations,
\begin{itemize}
\item a closed dynamic $\alpha:\mathbf{C}\rightharpoondown \mathbf{PY}$ has to satisfy those properties for objects $S$, $T$,... and composable arrows $e$, $d$ of $\mathbf{C}$:
\subitem [disjunctivity] $S\neq T \Rightarrow S^\alpha \cap T^\alpha=\emptyset$,
\subitem [lax composition] $(e\circ d)^\alpha \subseteq e^\alpha \odot d^\alpha$,
\subitem [lax identity] $(Id_S)^\alpha \subseteq Id_{S^\alpha}$,
\item a dynamorphism $\delta:\alpha\looparrowright \beta$ is such that 
\subitem [lax naturality] $\forall (S\stackrel{d}{\rightarrow}T)\in\overrightarrow{\mathbf{C}}, \delta_T\odot d^\alpha\subseteq d^\beta\odot \delta_S$.
\end{itemize}
\end{rmq}

A closed dynamic $\alpha:\mathbf{C}\rightharpoondown \mathbf{PY}$ is said to be \emph{deterministic (resp. quasi-deterministic)} if for every $(S\stackrel{d}{\rightarrow}T)\in\overrightarrow{\mathbf{C}}$ the transition $d^\alpha$ is deterministic (resp. quasi-deterministic).
 Likewise, a dynamorphism $\delta$ is said to be \emph{(quasi-)deterministic} if all transitions $\delta_S$ are (quasi-)deterministic.
 
\begin{df}
A deterministic closed dynamic is called a \emph{clock}.
\end{df}

\subsection{The category $\mathcal{CD}$ of closed dynamics}

\begin{df}
We define the category \emph{$\mathcal{CD}$ of all closed dynamics} taking as \emph{dynamorphisms} from $\alpha:\mathbf{C}\rightharpoondown \mathbf{PY}$ to $\beta:\mathbf{D}\rightharpoondown \mathbf{PY}$  couples $(\Delta,\delta)$ with $\Delta:\mathbf{C}\rightarrow\mathbf{D}$ a functor, and $\delta:\alpha \looparrowright (\beta\circ \Delta)$ a lax-natural transformation. 
\end{df}

\begin{rmq}
When $\mathbf{C}=\mathbf{D}$ we assume by default that $\Delta=Id_\mathbf{C}$. In other words, $\mathcal{CD}_\mathbf{C}$ is not a full subcategory of $\mathcal{CD}$.
\end{rmq}

\subsection{The category $\mathcal{MD}_{(\mathbf{C}, L)}$ of $L$ multi-dynamics on $\mathbf{C}$}

Let $L$ be a non-empty set. 

\begin{df}[$L$-multi-dynamics on $\mathbf{C}$] We denote by $\mathcal{MD}_{(\mathbf{C}, L)}$ the category 
\begin{itemize}
\item with objects, called  \emph{$L$-multi-dynamics on $\mathbf{C}$}, the disjunctive%
\footnote{
As in the closed dynamics case, a functor $\alpha:\mathbf{C}\rightharpoondown \mathbf{PY}^{\underrightarrow{\scriptstyle{L}}}$ is said to be disjunctive if for all objects $S$ and $T$ in $\mathbf{C}$,
$
S\neq T \Rightarrow S^\alpha \cap T^\alpha=\emptyset.
$
} 
lax functors 
\[
\alpha:\mathbf{C}\rightharpoondown 
\mathbf{PY}^{\underrightarrow{\scriptstyle{L}}},
\]
\item with arrows  $\delta:\alpha\looparrowright\beta$ --- called \emph{$({\mathbf{C}, L})$-dynamorphisms} --- given by\footnote{About the intersection used here, see below the remark \ref{rmq natural intersection lambda}.}
\[
\mathcal{MD}_{(\mathbf{C}, L)}(\alpha,\beta)
=
\bigcap_{\lambda\in L} 
\mathcal{CD}_\mathbf{C}(\alpha_\lambda,\beta_\lambda),
\]
where for each $\lambda\in L$ we denote $\alpha_\lambda$ the closed dynamic associated in an obvious way with the multi-dynamic $\alpha$ for this $\lambda$.
\end{itemize}
\end{df}

Thus, by definition, and with the notation introduced in section §\ref{subsection PYL}, an $L$-multi-dynamic associates with each duration $(S\stackrel{d}{\rightarrow}T)\in\overrightarrow{\mathbf{C}}$ an $L$-family $S^\alpha {\stackrel{d^\alpha}{\twosquig}}_L T^\alpha$.\\

\begin{rmq}\label{rmq natural intersection lambda} [About the expression $\bigcap_{\lambda\in L} 
\mathcal{CD}_\mathbf{C}(\alpha_\lambda,\beta_\lambda)$]
If $\mathbf{C}$ and $\mathbf{E}$ are categories, and $\alpha$, $\alpha'$, $\beta$ and $\beta'$ are functors $\mathbf{C}\rightarrow\mathbf{E}$ such that for all objects $C\in\dot{\mathbf{C}}$ we have $\alpha(C)=\alpha'(C)$ and $\beta(C)=\beta'(C)$, then if $(\alpha,\beta)\neq (\alpha',\beta')$, the sets of natural transformations  $\mathit{Nat}(\alpha,\beta)$ and $\mathit{Nat}(\alpha',\beta')$ are disjoint
\[
\mathit{Nat}(\alpha,\beta)\cap \mathit{Nat}(\alpha',\beta')=\emptyset
\]
because a natural transformation $\delta\in \mathit{Nat}(\alpha,\beta)$ has a unique domain $\alpha$ and a unique codomain $\beta$. Nevertheless, forgetting  this information about domains and codomains  and identifying $\delta$ just with a family of arrows $\left(\alpha(C)\stackrel{\delta_C}{\rightarrow}\beta(C)\right)_{C\in \dot{\mathbf{C}}}$ belonging to $\overrightarrow{\mathbf{E}}$, we'll write $\delta\in \mathit{Nat}(\alpha,\beta)$  to say only that $\delta$ is natural in respect with $\alpha$ and $\beta$, that is
\[
\forall(S\stackrel{d}{\rightarrow}T)\in\overrightarrow{\mathbf{C}}, \delta_T\circ \alpha(d)=\beta(d)\circ\delta_S.
\]
Then, $\mathit{Nat}(\alpha,\beta)\cap \mathit{Nat}(\alpha',\beta')$ will designate the set of families of $\mathbf{E}$-arrows $\delta=((\delta_C)_{C\in \dot{\mathbf{C}}})$ which are natural in respect together with $(\alpha,\beta)$ and with $(\alpha',\beta')$. In general, this kind of intersection will be non-empty.

The same remark can be made about lax transformations between lax functors between bicategories, and it is in this sense that we use expressions like $\bigcap_{\lambda\in L} 
\mathcal{CD}_\mathbf{C}(\alpha_\lambda,\beta_\lambda)$. To emphasize the forgetfulness of information about domains and codomains, we could use some notation like $\underline{\mathcal{CD}_\mathbf{C}}$ but it would be a little heavy, so we prefer not to do so.
\end{rmq}

\begin{rmq}
Unsurprisingly, an $L$-multi-dynamic $\alpha:\mathbf{C}\rightharpoondown \mathbf{PY}^{\underrightarrow{\scriptstyle{L}}}$ is said to be \emph{{(quasi-)}}\,\emph{deterministic}  if for every $\lambda\in L$ the closed dynamic $\alpha_\lambda$ is {(quasi-)}deterministic.
\end{rmq}

\subsection{The category $\mathcal{MD}_{\mathbf{C}}$ of multi-dynamics on $\mathbf{C}$} 

Let $\alpha:\mathbf{C}\rightharpoondown 
\mathbf{PY}^{\underrightarrow{\scriptstyle{L}}}$ 
and 
$\beta:\mathbf{C}\rightharpoondown 
\mathbf{PY}^{\underrightarrow{\scriptstyle{M}}}$
be  multi-dynamics on $\mathbf{C}$. 
We define a \emph{$\mathbf{C}$-dynamorphism} $\delta:\alpha\looparrowright\beta$ 
as a couple $(\theta,\underline{\delta})$ 
with $\theta:L\rightarrow M$ a map, 
and $\underline{\delta}\in\bigcap_{\lambda\in L} 
\mathcal{CD}_\mathbf{C}(\alpha_\lambda,\beta_{\theta(\lambda)})$. The lax natural part $\underline{\delta}$ of such a dynamorphism $\delta$ is often itself simply denoted $\delta$. Thus, to be a dynamorphism, $(\theta,\delta)$ must satisfy the lax naturality condition
\[\forall \lambda\in L,\forall (S\stackrel{d}{\rightarrow}T)\in\overrightarrow{\mathbf{C}}, \delta_T\odot d^\alpha_\lambda\subseteq d^\beta_{\theta(\lambda)}\odot \delta_S.\]

\begin{df}
We define the category $\mathcal{MD}_{\mathbf{C}}$ taking as objects all  multi-dynamics on $\mathbf{C}$, and as arrows all $\mathbf{C}$-dynamorphisms between them.
\end{df}

Naturally,  a dynamorphism $\delta$ is said to be \emph{(quasi-)deterministic} if for every object $S\in\dot{\mathbf{C}}$, $\delta_S$ is (quasi-)deterministic.

\begin{rmq}
Closed dynamics on $\mathbf{C}$ can be seen as a full subcategory of  $\mathcal{MD}_{\mathbf{C}}$, and we can in particular consider dynamorphisms between closed dynamics on $\mathbf{C}$ and multi-dynamics on $\mathbf{C}$.
\end{rmq}

\subsection{The category $\mathcal{MD}$ of multi-dynamics} 

Let $\alpha:\mathbf{C}\rightharpoondown 
\mathbf{PY}^{\underrightarrow{\scriptstyle{L}}}$ 
and 
$\beta:\mathbf{D}\rightharpoondown 
\mathbf{PY}^{\underrightarrow{\scriptstyle{M}}}$
be  multi-dynamics with possibly different engines. 
We define a \emph{dynamorphism} $\alpha\looparrowright\beta$ 
as a triple $(\theta,\Delta,\delta)$ 
with 
\begin{itemize}
\item $\theta:L\rightarrow M$ a map, 
\item $\Delta:\mathbf{C}\rightarrow\mathbf{D}$ a functor,
\item $\delta\in\bigcap_{\lambda\in L}\mathcal{CD}_\mathbf{C}(\alpha_\lambda,\beta_{\theta(\lambda)}\circ\Delta)$.
\end{itemize}

The last condition is equivalent to
\[
(\Delta,\delta)\in\bigcap_{\lambda\in L}\mathcal{CD}(\alpha_\lambda,\beta_{\theta(\lambda)}),
\]
\emph{i.e.} the lax naturality condition
\[\forall \lambda\in L,\forall (S\stackrel{d}{\rightarrow}T)\in\overrightarrow{\mathbf{C}}, \delta_T\odot d^\alpha_\lambda\subseteq (\Delta d)^\beta_{\theta(\lambda)}\odot \delta_S\]
has to be satisfied.

\begin{df}
The category $\mathcal{MD}$ of multi-dynamics is defined taking as objects all multi-dynamics, and as arrows all dynamorphisms between them.
\end{df}

\section{The category $\mathcal{OD}$ of open dynamics}

\subsection{Definition}

\begin{df}\label{df open dynamics}  An \emph{open dynamic} $A$ with engine $\mathbf{C}$ is the data
\[A=\left((\alpha:\mathbf{C}\rightharpoondown \mathbf{PY}^{\underrightarrow{\scriptstyle{L}}}) \stackrel{\rho}{\looparrowright}  (\mathbf{h}:\mathbf{C}\rightarrow \mathbf{Sets})\right)\] of
\begin{itemize}
\item a non-empty set $L$,
\item a multi-dynamic $\alpha\in \mathcal{MD}_{(\mathbf{C}, L)}\subset \mathcal{MD}_{\mathbf{C}}$,
\item a clock $\mathbf{h}\in \mathcal{CD}_{\mathbf{C}}\subset \mathcal{MD}_{\mathbf{C}}$,
\item a deterministic dynamorphism 
$
{\rho}\in \mathcal{MD}_{\mathbf{C}}(\alpha,\mathbf{h})
$ called \emph{datation}.
\end{itemize}
\end{df}

According to the definitions given in \textbf{§\,2.4.2} of \cite{Dugowson:20150809} and \textbf{§\,1.2.2} of \cite{Dugowson:20160831}, a \emph{dynamorphism} from an open dynamic
\[A=\left((\alpha:\mathbf{C}\rightharpoondown \mathbf{PY}^{\underrightarrow{\scriptstyle{L}}}) \stackrel{\rho}{\looparrowright}  (\mathbf{h}:\mathbf{C}\rightarrow \mathbf{Sets})\right)\]
to an open dynamic
\[B=\left((\beta:\mathbf{D}\rightharpoondown \mathbf{PY}^{\underrightarrow{\scriptstyle{M}}}) \stackrel{\tau}{\looparrowright}  (\mathbf{k}:\mathbf{D}\rightarrow \mathbf{Sets})\right)\]
is a quadruplet $(\theta, \Delta, \delta, \varepsilon)$ with
\begin{itemize}
\item $(\theta,\Delta,\delta)\in\mathcal{MD}(\alpha,\beta)$,
\item $(\Delta,\varepsilon)\in\mathcal{CD}(\mathbf{h},\mathbf{k})$,
\item this lax synchronization condition satisfied:
\[
\forall S\in\dot{\mathbf{C}}, 
\tau_{\Delta_S}\odot\delta_S\subseteq\varepsilon_S\odot\rho_S.
\]
\end{itemize}

We'll denote by $\mathcal{OD}$ the category of all open dynamics, with dynamorphisms as arrows.

\subsection{Realizations of an open dynamic}

With the definitions given above, the other notions of our systemic theory of interactivity keep the same definitions as previously given in \cite{Dugowson:20160831}: parametric quotients, interactions, connectivity structures of an interaction, dynamical families, open dynamics generated by a dynamical family, ...

In particular, given an open dynamic $A=\left((\alpha:\mathbf{C}\rightharpoondown \mathbf{PY}^{\underrightarrow{\scriptstyle{L}}}) \stackrel{\rho}{\looparrowright}  (\mathbf{h}:\mathbf{C}\rightarrow \mathbf{Sets})\right)$ we have 

\begin{df}
A \emph{realization} (or a \emph{solution}) of $A$ is a quasi-deterministic dynamorphism $(\mathfrak{s}:\mathbf{h}\looparrowright\alpha)\in \mathcal{MD}_\mathbf{C}(\mathbf{h},\alpha)$ satisfying the lax condition $\rho\odot \mathfrak{s}\subseteq Id_\mathbf{h}$. 
\end{df}

In other words\footnote{See \cite{Dugowson:20160831}, \textbf{§\,1.3.1.}}, such a realization is a couple $\mathfrak{s}=(\lambda,\sigma)$ with $\lambda\in L$ and $\sigma\in \mathcal{CD}_\mathbf{C}(\mathbf{h},\alpha_\lambda)$ which is a partial function $\sigma:st(\mathbf{h})\dashrightarrow st(\alpha)$ defined on a subset $D_\sigma\subset st(\mathbf{h})$ and satisfying these properties:

\[\forall t\in D_\sigma, \rho(\sigma(t))=t,\]
\[\forall S\in\dot{\mathbf{C}},\forall t\in S^\mathbf{h}\cap  D_\sigma, \sigma(t)\in S^\alpha,\]
\[\forall (S\stackrel{d}{\rightarrow}T)\in\overrightarrow{\mathbf{C}}, 
 \forall t\in S^\mathbf{h}, 
\left( d^\mathbf{h}(t)\in D_\sigma\Rightarrow 
\left( t\in D_\sigma 
\,\mathrm{and}\,
\sigma(d^\mathbf{h}(t))\in d^\alpha_\lambda (\sigma(t))\right)\right).\]

\vspace{0.7cm}

\begin{center}
***
\end{center}

\vspace{0.5cm}

\textbf{Thanks} to Andrée Ehresmann and Mathieu Anel for having encouraged me to reword the previously developed notion of ``sub-functorial dynamic''  in terms of lax-functors.

\bibliographystyle{plain}




\tableofcontents
\end{document}